\documentclass{tac}
  \usepackage[all,2cell,arrow]{xy}
  \usepackage{amsfonts}
  \usepackage{amsmath}
  \usepackage{amssymb}
 \usepackage{color}
 \usepackage{graphicx}

\usepackage{enumerate,xspace}

\UseAllTwocells
\CompileMatrices

\renewcommand{\epsilon}{\varepsilon}
\renewcommand{\phi}{\varphi}

\newcommand{\ca}{\ensuremath{\mathcal A}\xspace}
\newcommand{\cb}{\ensuremath{\mathcal B}\xspace}
\newcommand{\cc}{\ensuremath{\mathcal C}\xspace}
\newcommand{\cd}{\ensuremath{\mathcal D}\xspace}
\newcommand{\ce}{\ensuremath{\mathcal E}\xspace}
\newcommand{\cf}{\ensuremath{\mathcal F}\xspace}

\newcommand{\ck}{\ensuremath{\mathcal K}\xspace}

\newcommand{\cm}{\ensuremath{\mathcal M}\xspace}

\newcommand{\cv}{\ensuremath{\mathcal V}\xspace}
\newcommand{\cw}{\ensuremath{\mathcal W}\xspace}
\newcommand{\cx}{\ensuremath{\mathcal X}\xspace}

\newcommand{\bba}{\ensuremath{\mathbb A}\xspace}
\newcommand{\bbb}{\ensuremath{\mathbb B}\xspace}

\newcommand{\bbd}{\ensuremath{\mathbb D}\xspace}

\newcommand{\bbf}{\ensuremath{\mathbb F}\xspace}

\newcommand{\bbk}{\ensuremath{\mathbb K}\xspace}

\newcommand{\bbo}{\ensuremath{\mathbb O}\xspace}

\newcommand{\bbs}{\ensuremath{\mathbb S}\xspace}

\newcommand{\bbz}{\ensuremath{\mathbb Z}\xspace}

\newcommand{\bbbo}{\ensuremath{\bbb\bbo}\xspace}

\newcommand{\enr}{\ensuremath{\textnormal{-}\Cat}\xspace}

\newcommand{\BO}{\ensuremath{\mathbf{BO}}\xspace}

\newcommand{\Cat}{\ensuremath{\mathbf{Cat}}\xspace}
\newcommand{\DblCat}{\ensuremath{\mathbf{DblCat}}\xspace}
\newcommand{\Set}{\ensuremath{\mathbf{Set}}\xspace}

\newcommand{\op}{\ensuremath{^{\textnormal{op}}}}

\DeclareMathOperator{\disc}{disc}
\DeclareMathOperator{\cod}{cod}
\DeclareMathOperator{\dom}{dom}
\DeclareMathOperator{\id}{id}
\DeclareMathOperator{\Id}{Id}
\DeclareMathOperator{\EM}{EM}

\DeclareMathOperator{\core}{core}
\DeclareMathOperator{\Fam}{Fam}

\DeclareMathOperator{\KL}{KL}
\DeclareMathOperator{\Mnd}{Mnd}

\newcommand{\two}{\ensuremath{\mathbf{2}\xspace}}

\begin{document}

 \title[The universal property of the 2-category of monads]{What is the universal property of the 2-category of monads?}

\author{Stephen Lack and Adrian Miranda}
\address{School of Mathematical and Physical Sciences, Macquarie University NSW 2109, Australia}
\eaddress{steve.lack@mq.edu.au\CR adrian.miranda@hdr.mq.edu.au}
\thanks{The first-named author acknowledges with gratitude the support of an
 Australian Research Council Discovery Project DP190102432, and the
 second-named author an MQRES PhD Scholarship, 20192497.\\
 We are happy to acknowledge useful conversations with John Power on
the topic of the paper; AM is also grateful to Charles Walker for
discussions around the connection with double categories.}
\dedication{In memory of our colleague Pieter Hofstra}
\copyrightyear{2022}


\maketitle


\begin{abstract}
For a 2-category $\ck$, we consider Street's 2-category $\Mnd(\ck)$ of
monads in $\ck$, along with Lack and Street's 2-category $\EM(\ck)$
and the identity-on-objects-and-1-cells 2-functor
$\Mnd(\ck)\to\EM(\ck)$ between them. We show that this 2-functor can
be obtained as a ``free completion'' of the 2-functor $1\colon
\ck\to\ck$. We do this by regarding 2-functors which act as the
identity on both objects and 1-cells as categories enriched a
cartesian closed category $\BO$ whose objects are identity-on-objects
functors. We also develop some of the theory of $\BO$-enriched
categories. 
\end{abstract}

\section{Introduction}

In \cite{ftm}, Street introduced, for a given 2-category $\ck$,  a
2-category $\Mnd(\ck)$ whose objects are the monads in $\ck$, and
showed how various aspects of the theory of monads can be understood
in terms of this construction. For example, there is a 2-functor
$\Id\colon\ck\to\Mnd(\ck)$ sending each object of $\ck$ to the identity monad
on that object, and this $\Id$ has a right adjoint just when $\ck$
admits what later came to be known as Eilenberg-Moore objects (which in
the classical case, $\ck=\Cat$, correspond to the category of algebras
for the monad, introduced by Eilenberg and Moore
\cite{EilenbergMoore-Triples}).

Some thirty years later, in \cite{ftm2}, a variant $\EM(\ck)$ of
$\Mnd(\ck)$ was introduced, having the same objects and 1-cells, but
a different notion of 2-cell. There is a 2-functor
$\Mnd(\ck)\to\EM(\ck)$ between them, which acts as the identity on
objects and 1-cells. As shown in \cite{ftm2}, it is still the
case that the composite 2-functor $\ck\to\Mnd(\ck)\to \EM(\ck)$ has a
right adjoint just when $\ck$ admits Eilenberg-Moore objects, but now
there is a conceptual reason: $\EM(\ck)$ is the free completion of
$\ck$ under Eilenberg-Moore objects.

One might then ask whether the original $\Mnd(\ck)$ itself has a
universal property, and lo these twenty years after \cite{ftm2} we offer our
response to this question.

In fact, as often happens in mathematics, rather than answer this question
as is, we first reformulate it, and then respond to the modified
question. Thus rather than seek to characterize $\Mnd(\ck)$
alone in terms of some universal property, we consider $\Mnd(\ck)$ and
$\EM(\ck)$ together, along with the 2-functor $\Mnd(\ck)\to\EM(\ck)$
which connects them, which as mentioned above {\em acts as the identity on both
objects and 1-cells}. We consider 2-functors with this latter property as a structure in
its own right, and describe a universal property of
$\Mnd(\ck)\to\EM(\ck)$ relative to the identity 2-functor $\ck\to\ck$.

The way we do this is to consider such 2-functors as a certain sort of
enriched category. In more detail, there is a cartesian closed
category $\BO$ whose objects are identity-on-objects functors and
whose morphisms are commutative squares of functors. A $\BO$-enriched
category is essentially the same as a 2-functor which acts as the identity
on objects and on 1-cells. 

We then develop a little of the theory of $\BO$-enriched category
theory, including in particular weighted limits and colimits, and free
completions under these. We also show how every $\Cat$-enriched
weight gives rise to a corresponding $\BO$-enriched weight (actually
in two different ways), and so in particular we have a $\BO$-enriched
notion of Eilenberg-Moore object.

Our answer to the (reformulated) problem, then, is that
$\Mnd(\ck)\to\EM(\ck)$ is the free completion of $1\colon\ck\to\ck$
under these $\BO$-enriched Eilenberg-Moore objects.

As was the case in \cite{ftm2}, it is technically more convenient to
deal with colimits rather than limits when it comes to free
completions, so we actually work with the colimit notion corresponding
to Eilenberg-Moore objects, namely Kleisli objects.

We begin, in Section~\ref{sect:BO}, by introducing our cartesian
closed category $\BO$, and describing various relationships it has to
$\Cat$. Then in Section~\ref{sect:enrichment} we begin our study of
$\BO$-enriched categories, including enriched presheaf categories, as
well as various ``change-of-base'' constructions linking
$\BO$-categories with other sorts of enriched categories. In
Section~\ref{sect:colimits}, we study various examples of
$\BO$-enriched colimits, leading up to our main result, Theorem~\ref{thm:EM}, characterizing
$\Mnd(\ck)\to\EM(\ck)$ as a free completion.


\section{The cartesian closed category $\BO$}\label{sect:BO}

We write $\Cat_1$ for the cartesian closed category of (small)
categories and functors, where the subscript $1$ is to emphasize that
we are thinking of this as a mere category rather than a 2-category.
The arrow category $\Cat^\two_1$ is also cartesian closed. We write $\BO$
for the full subcategory of $\Cat^\two_1$ consisting of the
identity-on-object functors.\footnote{
  Very little would change if we were to work with functors which are
  merely bijective on objects, and indeed these are the source of our
  notation $\BO$. Of course every bijective-on-object functor is
  isomorphic in $\Cat^\two_1$ to an object of $\BO$.}
Our standard notation for objects of $\BO$ is a letter like $A$, then
we write $e_A\colon A_t\to A_\ell$ for the corresponding
functor\footnote{The subscripts $t$ and $\ell$ stand for tight and
  loose, as in \cite{enhanced}; see also
  Section~\ref{sect:BO-enriched} below.},
sometimes dropping the subscript $A$ in $e_A$. A typical morphism
$f\colon A\to B$ has the form of a commutative square
\[ \xymatrix{
    A_t \ar[r]^-{f_t} \ar[d]_{e_A} & B_t \ar[d]^{e_B} \\ A_\ell
    \ar[r]_-{f_\ell} & B_\ell. } \]
Given any functor $g\colon X\to Y$ there is an essentially unique
factorization (often called the {\em bo-ff factorization})
\[ \xymatrix{
    X \ar[r]^-{e} & Z\ar[r]^-{j} & Y } \]
with $e$ the identity on objects and $j$ fully faithful.

\begin{proposition}
  The full subcategory $\BO$ of $\Cat^\two_1$ is both reflective and coreflective.
\end{proposition}

\begin{proof}
  The {\em coreflection} of an object $g\colon X\to Y$ is given by the
  identity-on-objects part $e\colon X\to Z$ of the bo-ff
  factorization; indeed, if $(\ce,\cm)$ is a
  factorization system on a category $\cc$, then $\ce$ determines a
  full coreflective subcategory of $\cc^\two$ and the coreflection
  sends a morphism to the $\ce$-part of the factorization.

  The {\em reflection} of $g\colon X\to Y$ is given by the induced map
  $g'\colon X\to Y$ in the diagram 
\[ \xymatrix{
    X_0 \ar[r]^-e \ar[d]_{g_0} & X \ar[d]_{\eta} \ar@/^1pc/[ddr]^g \\
    Y_0 \ar[r] \ar@/_1pc/[drr]_-e  & X' \ar[dr]^{g'} \\
    && Y } \]
in which the arrows labelled $e$ are identity-on-object inclusions of
discrete subcategories, and the square is a pushout.  
\end{proof}

\begin{remark}\label{rmk:2-category}
  Of course $\Cat_1^\two$ can be made into a 2-category $\Cat^\two$,
  and so $\BO$ becomes a full sub-2-category. As such, it is still coreflective, but not
  reflective; indeed it is not closed under powers.
\end{remark}

It follows from the coreflectivity and the fact that $\BO$ is
closed in $\Cat^\two_1$ under finite products (which in turn follows
from reflectivity or can easily be checked) that
$\BO$ is, like $\Cat_1^\two$, cartesian closed, with the internal hom in
$\BO$ formed as the coreflection of the internal hom in $\Cat_1^\two$.

Explicitly, if $e_A\colon A_t\to A_\ell$ and $e_B\colon
B_t\to B_\ell$ are in \BO, the internal hom $[A,B]$ has:
\begin{itemize}
\item objects are commutative squares (in other words, morphisms
  $f=(f_t,f_\ell)\colon A\to B$ in $\BO$);
\item $[A,B]_t(f,g)$ consists of natural transformations $f_t\to
  g_t$ and $f_\ell\to g_\ell$ subject to the obvious compatibility
  condition (these are also the 2-cells if $\BO$ is made into a
  2-category, as in Remark~\ref{rmk:2-category}).
\item $[A,B]_\ell(f,g)$ consists just of natural transformations $f_\ell\to g_\ell$.
\end{itemize}

As well as the adjunctions between $\BO$ and $\Cat_1^\two$, there are
also adjunctions between $\BO$ and $\Cat_1$. In particular, the codomain
functor $\cod\colon \BO\to \Cat_1$ has a left adjoint $\disc\colon\Cat_1\to\BO$ sending a
category $C$ to the inclusion $C_0\to C$ of the discrete category with
the same objects as $C$. Once again, this is just an adjunction of
ordinary categories, not of 2-categories. On the other hand, both
adjoints preserve finite products and so this is a monoidal
adjunction, and induces a 2-adjunction $\disc_*\dashv\cod_*$ between $\BO$-categories and 2-categories.

In fact there is a chain of monoidal 
adjunctions\footnote{An adjunction between cartesian
  closed categories is monoidal just when the left adjoint preserves
  finite products.}
\[\pi\dashv\disc\dashv\cod\dashv\id\dashv\dom\dashv\textnormal{ch} \]
between $\BO$ and $\Cat_1$. It is $\disc\dashv\cod$ and
$\id\dashv\dom$ which will be particularly important in what
follows. Here $\dom\colon\BO\to\Cat_1$ sends an identity-on-objects functor to its
domain, while $\id\colon\Cat_1\to\BO$ sends a category to the corresponding identity functor.

\begin{remark}\label{rmk:pi0-fpp}
  The adjunction $\pi\dashv \disc$ will not be significant in what
  follows, so we shall not describe it in detail. It sends an
  identity-on-objects map $A\to A'$ to its pushout with the canonical
  map $A\to \pi_0A$.
In order to see that this $\pi$ preserves finite products, and so that 
$\pi\dashv\disc$ is monoidal, it is perhaps
easiest to use the Day reflection theorem \cite[Theorem~1.2 and
Corollary~2.1]{Day-reflection}, according to which it will suffice to check that if $A\in\BO$ and
$C\in\Cat_1$, then the internal hom $[A,\disc C]$ is in the image of
$\disc$. To do this, one observes that an object $B$ of $\BO$ is in the image of
$\disc$ if and only if $B_t$ is discrete, and then uses the explicit description of $[A,\disc C]_t$ given above. 
\end{remark}

\begin{remark}
  We have directly described a cartesian closed structure on $\BO$ and
  in the following section we shall consider enrichment over this
  structure. But it is also worth pointing out that $\BO$ can itself
  be understood in terms of enrichment. The category $\Set^\two$ is
  also cartesian closed, and the category of $\Set^\two$-enriched
  categories and $\Set^\two$-enriched functors is isomorphic to $\BO$:
  see \cite[Example~1]{Power-Premonoidal}. This immediately implies
  that $\BO$ is cartesian closed; furthermore, the adjunctions
  $\cod\dashv\id\dashv\dom\dashv\textnormal{ch}$ between $\BO$ and $\Cat_1$ can
  themselves be seen as arising via change of base-of-enrichment from
  monoidal adjunctions between $\Set^\two$ and $\Set$.
\end{remark}


\begin{remark}
  Another point of view on $\BO$ relates it to the cartesian closed
  category $\DblCat$ of double categories. There is a functor
  $\BO\to\DblCat$ sending an identity on objects functor to its
  “higher kernel”; this functor is fully faithful, with image given by
  the cateads, and moreover has a finite-product-preserving
  reflection. See \cite{BournPenon} for the notion of cateads and
  their equivalence to bijective on objects functors in a more general
  context, or \cite{Bourke-thesis} for further details, including a simpler
  treatment for the specific case at hand. The fact that the
  reflection preserves finite products can be found in
  \cite[Proposition~8.30]{Bourke-thesis}. 
\end{remark}

\section{$\BO$-enriched categories}\label{sect:enrichment}

Since $\BO$ is cartesian closed, we can consider enrichment over it.
In this section we investigate what some of the basic theory of
enriched categories looks like in the case of enrichment over $\BO$. 

\subsection{$\BO$-enriched categories}\label{sect:BO-enriched}

We generally use blackboard bold for the names of $\BO$-categories. 
A $\BO$-category $\bba$ has objects $A,B,C$, and so on; and, for objects
$A$ and $B$, a $\BO$-valued hom $\bba(A,B)$. This consists of an
identity-on-objects functor, which we write as 
$E_\bba\colon\ca_t(A,B)\to \ca_\ell(A,B)$.
The composition and identities make $\ca_t(A,B)$ and $\ca_\ell(A,B)$ into the
hom-categories of 2-categories $\ca_t$ and $\ca_\ell$, in such a way
that the $E_\bba$ define a 2-functor $E_\bba\colon\ca_t\to\ca_\ell$
which acts as the identity on objects and on 1-cells.  We often drop the subscript
$\bba$ and write simply $E$. Conversely, any 2-functor which acts as
the identity on objects and on 1-cells arises in this way from a
unique $\BO$-category. We shall routinely identify
$\BO$-categories with the corresponding 2-functors.

Following the naming convention of \cite{enhanced}, 2-cells of
$\ca_t$ are referred to as {\em tight}, and 2-cells of $\ca_\ell$ as
{\em loose}; this is also the origin of the subscripts $t$ and $\ell$.

Given $\BO$-categories $\bba$ and $\bbb$, seen as 2-functors
$E_\bba\colon \ca_t\to\ca_\ell$ and $E_\bbb\colon \cb_t\to\cb_\ell$, a
$\BO$-functor $F\colon\bba\to\bbb$ consists of 2-functors
$F_t\colon\ca_t\to\cb_t$ and $F_\ell\colon\ca_\ell\to\cb_\ell$ making
the square
\[ \xymatrix{
    \ca_t \ar[r]^-{F_t} \ar[d]_{E_\bba} & \cb_t \ar[d]^{E_\bbb} \\
    \ca_\ell \ar[r]_-{F_\ell} & \cb_\ell } \]
commute.

Similarly, if $G\colon\bba\to\bbb$ is also a $\BO$-functor, a
$\BO$-natural transformation $F\to G$ consists of a pair of 2-natural
transformations $\phi_t\colon F_t\to G_t$ and $\phi_\ell\colon
F_\ell\to G_\ell$ satisfying the evident compatibility condition
$E_\bbb \phi_t = \phi_\ell E_\bba$.

We can summarize the results of this analysis using the 2-category $2\enr$ of
2-categories, 2-functors, and 2-natural transformations, and the
2-category $2\enr^\two$ of arrows in $2\enr$.

\begin{proposition}
  The 2-category $\BO\enr$ of $\BO$-categories, $\BO$-functors, and
  $\BO$-natural transformations is isomorphic to the full
  sub-2-category of $2\enr^\two$ consisting of those 2-functors which
  act as the identity on objects and on 1-cells. 
\end{proposition}

Every enriched category has an underlying ordinary category; the
underlying ordinary category of a $\BO$-category $\bba$ is the
underlying ordinary category of the 2-category $\ca_t$, which is
of course the same as the underlying ordinary category of the 2-category $\ca_\ell$.

\subsection{The $\BO$-category $\bbbo$.}

We write $\bbbo$ for the $\BO$-category coming from the cartesian
closed structure of $\BO$ itself. The corresponding 2-functor
$\BO_t\to\BO_\ell$ can be described as follows.

The objects of $\BO_t$, $\BO_\ell$, $\bbbo$, and $\BO$ all coincide,
and are just the identity-on-objects functors. Once again, the 1-cells
of $\BO_t$, $\BO_\ell$, and $\BO$ all coincide (and we could add
$\bbbo$ to this list if we understand the 1-cells of a $\BO$-category
to be the 1-cells of the underlying ordinary category): these are all
just the  commutative
squares of functors, with vertical maps acting as the identity on
objects.

Given 1-cells $F=(F_t,F_\ell)$ and $G=(G_t,G_\ell)$ from
$E_\bba\colon \ca_t\to\ca_\ell$ to $E_\bbb\colon\cb_t\to\cb_\ell$, a
2-cell $F\to G$ in $\BO_t$ consists of natural transformations
$\phi_t\colon F_t\to G_t$ and $\phi_\ell\colon F_\ell\to G_\ell$
satisfying the evident compatibility condition $\phi_\ell E_\bba=
E_\bbb \phi_t$. Thus in fact $\BO_t$ can be identified with the full
sub-2-category of $\Cat^\two$ consisting of those functors which act
as the identity on objects.

A 2-cell $F\to G$ in $\BO_\ell$ consists of just the single natural
transformation $\phi_\ell\colon F_\ell\to G_\ell$, and $E_{\bbbo}$
sends a 2-cell $(\phi_t,\phi_\ell)$ to $\phi_\ell$. Thus in fact we
have a commutative square
\[ \xymatrix{
    \BO_t \ar[r]^-{H} \ar[d]_{E_\bbbo} & \Cat^\two \ar[d]^{\cod} \\
    \BO_\ell \ar[r]_-{\cod'}  & \Cat } \]
where the upper horizontal $H$ is the fully faithful inclusion mentioned
above, the left vertical acts as the identity on objects and 1-cells,
and the lower horizontal $\cod'$ is fully faithful on 2-cells. This is enough
to determine $\BO_\ell$.

\subsection{Presheaves}

A presheaf on $\bba$ consists of a $\BO$-functor $\bba\op\to\bbbo$,
but we can also analyze this in terms of the 2-functor $E_\bba\colon
\ca_t\to\ca_\ell$. Given such a presheaf, we obtain a commutative
diagram
\[ \xymatrix{
    \ca\op_t \ar[d]_{E} \ar[r] & \BO_t \ar[r]^-{H} \ar[d]_{E_{\bbbo}} &
    \Cat^\two \ar[d]^{\cod} \\
    \ca\op_\ell \ar[r] & \BO_\ell \ar[r]_-{\cod'} & \Cat } \]
of 2-functors, and by the universal property of the power $\Cat^\two$,
this in turn determines a diagram
\[ \xymatrix @R1pc {
    \ca\op_t \ar[dd]_E \ar[dr]^{P_t}_{~}="1" \\
    & \Cat \\
    \ca\op_\ell \ar[ur]_{P_\ell}^{~}="2"
    \ar@{=>}"1";"2"^{P_E} } \]
where $P_\ell$ is the composite lower horizontal in the previous
diagram, while the composite upper horizontal $\ca\op_t\to\Cat^\two$
corresponds to the 2-natural map $P_E\colon P_t\to P_\ell E$. 

When does a 2-natural $P_E\colon P_t\to P_\ell E$ arise in this way
from some $\BO$-presheaf? Such a $P_E$ does determine a unique map
$\overline{P}\colon\ca\op_t\to\Cat^\two$ with $\cod\overline{P}=P_\ell
E$. This $\overline{P}$
will land in $\BO_t$ just when the components of $P_E$ act as the
identity on objects. And indeed when this happens, we have the solid
part of the diagram 
\[ \xymatrix{
    \ca\op_t \ar[d]_{E} \ar[r] & \BO_t \ar[r]^-{E_{\bbbo}} & \BO_\ell \ar[d]^{\cod'} \\
    \ca\op_\ell \ar[rr] \ar@{.>}[urr] && \Cat} \]
in which the left vertical is the identity on objects and 1-cells, and the right
vertical is fully faithful on 2-cells. It follows that there is a unique
induced diagonal filler.

This proves:

\begin{proposition}\label{prop:presheaves}
Let $\bba$ be a $\BO$-category and $E\colon\ca_t\to\ca_\ell$ the
corresponding 2-functor. To give a presheaf on $\bba$ is equivalently
to give $\Cat$-presheaves $P_t$ and $P_\ell$ on $\ca_t$ and $\ca_\ell$
respectively, along with a 2-natural transformation $P_E\colon P_t\to
P_t E$ whose components are identity-on-objects functors. \endproof
\end{proposition}

From this point of view, a representable presheaf $\bba(-,A)$ has the
form
\[ \xymatrix @R1pc {
    \ca\op_t \ar[dd]_{E} \ar[dr]^{\ca_t(-,A)}_{~}="1" \\
    & \Cat. \\
    \ca\op_\ell \ar[ur]_{\ca_\ell(-,A)}^{~}="2"
    \ar@{=>}"1";"2"^{E} } \]

\subsection{Enriched presheaf categories}\label{sect:presheaf-cats}

Since the cartesian closed category $\BO$ is complete, we can construct enriched presheaf categories as in \cite[Section~2.2]{Kelly-book}.

Given two presheaves $P$ and $Q$ on the $\BO$-category $\bba$, there
is a \BO-valued hom $[\bba\op,\bbbo](P,Q)$. An object is a $\BO$-natural $f\colon P\to Q$:
this amounts to maps $f_t\colon P_t\to Q_t$ and $f_\ell\colon P_\ell\to
Q_\ell$ making the square
\[ \xymatrix{
P_t \ar[r]^-{f_t} \ar[d]_{P_E} & Q_t \ar[d]^{Q_j} \\
P_\ell J \ar[r]_-{f_\ell J} & Q_\ell J } \]
commute. A tight map $f\to g$ consists of compatible 2-cells $f_t\to
g_t$ and $f_\ell\to g_\ell$, while a loose map consists of just a
2-cell $f_\ell\to g_\ell$.

A more abstract way to summarize this is as follows. First form the pullback of
2-categories
\[ \xymatrix @C3pc {
    [\bba\op,\bbbo]_t \ar[r] \ar[d] & [\ca\op_\ell,\Cat] \ar[d]^{P^*_E}
    \\
    [\ca\op_t,\BO_t] \ar[r]_-{[\ca\op_t,\cod]} & [\ca\op_t,\Cat] } \]
and now factorize the upper horizontal as 
\[ \xymatrix @C3pc {
    [\bba\op,\bbbo]_t \ar[r]^-{E_{[\bba\op,\bbbo]}} & [\bba\op,\bbbo]_\ell \ar[r] &
    [\ca\op_\ell,\Cat] } \]
where $E_{[\bba\op,\bbbo]}$ acts as the identity on objects and on 1-cells,
and the other map is fully faithful on 2-cells.

\subsection{Change of base}

As observed in Section~\ref{sect:BO} above, there is a chain of  monoidal adjunctions
$\disc\dashv\cod\dashv\id\dashv\dom$ between $\BO$ and $\Cat_1$. These
induce 2-adjunctions $\disc_*\dashv\cod_*\dashv\id_*\dashv\dom_*$
between $\BO\enr$ and $2\enr$.

Given a $\BO$-category $\bba$, the 2-category $\dom_*\bba$ is $\ca_t$,
while $\cod_*\bba$ is $\ca_\ell$.

Given a 2-category $\cd$, the
$\BO$-category $\id_*\cd$ is the one corresponding to the identity
2-functor $1\colon\cd\to\cd$. We sometimes identify a 2-category
$\cd$ with the corresponding $\BO$-category $\id_*\cd$, and call a
$\BO$-category of this form {\em tight}, since it has the property
that every loose 2-cell has a unique tight structure (thus we might
``loosely'' say that all 2-cells are tight!)

On the other hand, $\disc_*\cd$ is the $\BO$-category corresponding to
the 2-functor $\cd_1\to\cd$, where $\cd_1$ is the underlying ordinary
category of $\cd$, seen as a locally discrete 2-category (no
non-identity 2-cells). We call a $\BO$-category of this form {\em
  loose}, since there are no non-identity tight 2-cells.

The tight and loose $\BO$-categories are analogous, respectively, to
the {\em chordate} and {\em inchordate} $\cf$-categories of
\cite{enhanced}; see also \cite[Example~16.2]{GarnerShulman-FreeCocompletion}.

\subsection{Change of base and tight presheaves}\label{sect:CoBpresh1}

As well as these 2-adjunctions between $\BO\enr$ and $2\enr$, there
are also various connections between the enriched presheaf categories.

First observe that the monoidal adjunction $\id\dashv\dom$ also induces a
2-adjunction between $\BO_t=\dom_*\bbbo$ and $\Cat$: recall that
$\BO_t$ is just the full sub-2-category of $\Cat^\two$ consisting of
the identity-on-objects functors.
For any 2-category $\cd$, we have
\[ \dom_*[\id_*\cd\op,\bbbo] \cong [\id_*\cd\op,\bbbo]_t \cong
  [\cd\op,\BO_t] \]
and so $\id\dashv\dom$ induces a 2-adjunction between $[\cd\op,\Cat]$ and
this $\dom_*[\id_*\cd\op,\bbbo]$. In particular, a 2-functor
$F\colon\cd\op\to\Cat$ will be sent to the $\BO$-functor
$F^{\id}\colon\id_*\cd\op\to\bbbo$ with $F^{\id}_t=F^{\id}_\ell=F$ and
$F_E\colon F_t\to F_\ell$ equal to the identity.

\begin{definition} A weight of the form $F^{\id}\colon \id_*\cd\op\to\bbbo$ for a
$\Cat$-weight $F\colon\cd\op\to\Cat$ will be called a {\em tight
  weight}. \label{defn:tight-weight}
\end{definition}

\begin{remark} \label{rmk:cob-presheaf}
  In fact these tight weights arise by general change of
  base arguments. If $\phi\colon\cv\to\cw$ is a monoidal functor, then
  any $\cv$-weight $F\colon\cd\op\to\cv$ determines a $\cw$-weight as
  follows. First apply change of base along $\phi$ to obtain a
  $\cw$-functor   $\phi_*F\colon \phi_*\cd\op\to\phi_*\cv$; now
  compose with the $\cw$-functor $\phi\colon \phi_*\cv\to\cw$. If
  moreover $\phi$ has a monoidal right adjoint $\psi$, then this process 
  defines a $\cv$-functor $[\cd\op,\cv]\to \psi_*[\phi_*\cd\op,\cw]$
  which in turn has a right $\cv$-functor induced by $\psi$. Applying
  this in the case of the monoidal adjunction $\id\dashv\dom$ sends a
  $\Cat$-weight $F$ to the corresponding tight $\BO$-weight $F^{\id}$.
\end{remark}

\subsection{Change of base and loose presheaves}\label{sect:CoBpresh2}

This time we consider the monoidal adjunction $\disc\dashv\cod$.  Once
again, we start with a $\Cat$-weight $F\colon\cd\op\to\Cat$. This time
we apply Remark~\ref{rmk:cob-presheaf} using $\disc\dashv\cod$ to
obtain a $\BO$-weight $\bbf\colon\disc_*\cd\op\to\bbbo$.

Recall that $\disc_*\cd$ corresponds to the inclusion 2-functor
$E\colon\cd_1\to\cd$, where $\cd_1$ is obtained from $\cd$ by
discarding all non-identity 2-cells. Then $F_t\colon\cd\op_1\to\Cat$
is defined to send $D\in\cd_1$ to the discrete category with the same
set of objects as $FD$. The identity-on-objects inclusions
$F_tD\to FD$ are the components of a 2-natural transformation
$F_E\colon F_t\to FE$ whose components act as the identity on
objects. This in turn determines the presheaf $\bbf$.

\begin{definition} A weight of the form
  $\bbf\colon \disc_*\cd\op\to\bbbo$ for a $\Cat$-weight
  $F\colon\cd\op\to\Cat$ will be called a {\em loose weight}. \label{defn:loose-weight}
\end{definition}

\subsection{Another example of change-of-base}

In this section we include, for interest's sake, a further example of
change-of-base. It will be not be used in the remainder of the paper.

If $\ck$ is a 2-category then we obtain another 2-category $\ck_g$
with the same objects and morphisms by discarding all non-invertible
2-cells. Thus $\ck_g\to\ck$ can be seen as a $\BO$-category.

In fact this also arises through a change-of-base process. There is a functor
$\core\colon\Cat_1\to\BO$ sending a category $A$ to the inclusion
$A_g\to A$, where $A_g$ is the subcategory of $A$ consisting of the
isomorphisms. This has a left adjoint $\pi$ sending $B_t\to B_\ell$ to the pushout
\[ \xymatrix{
    B_t \ar[r] \ar[d] & \pi_1(B_t) \ar[d] \\
    B_\ell \ar[r] & \overline{B} } \]
where $B_t\to \pi_1(B_t)$ is the map which universally inverts all
morphisms in $B_t$. The counit of $\pi\dashv\core$ is invertible.

Much as in Remark~\ref{rmk:pi0-fpp}, the Day reflection theorem \cite[Theorem~1.2 and
Corollary~2.1]{Day-reflection} can be used to see that $\pi$ preserves
finite products and so $\pi\dashv\core$ is monoidal, by observing that
$[B,\core A]$ is in the image of $\core$ for any $B\in\BO$ and
$A\in\Cat_1$.

\section{$\BO$-enriched colimits}\label{sect:colimits}

Let $M\colon\bba\op\to\bbbo$ be a presheaf, corresponding as in
Proposition~\ref{prop:presheaves} to $M_t\colon\ca\op_t\to\Cat$,
$M_\ell\colon\ca\op_\ell\to\Cat$, and
$M_E\colon M_t\to M_\ell E_\bba$. By the Yoneda lemma, for an object
$A\in\bba$ there is a bijection between $\BO$-natural
$\widehat{a}\colon\bba(-,A)\to M$ and maps $a\colon 1\to MA$ in
$\bbbo$, or equivalently of $a\in M_tA$ (or $a\in M_\ell A$). Then $M$
is representable when there is an $A\in\bba$ and $a\in M_tA$ as
above, for which the induced $\widehat{a}$ is invertible.  This in
turn is equivalent to invertibility of the 2-natural transformations
$\widehat{a}_t\colon \ca_t(-,A)\to M_t$ and
$\widehat{a}_\ell\colon\ca_\ell(-,A)\to M_\ell$. We call the
invertibility of $\widehat{a}_t$ and $\widehat{a}_\ell$, respectively,
{\em the tight and the loose aspects of the universal property}.

Recall \cite[Section~3.1]{Kelly-book} that if $\bbf\colon\bbd\to\bbbo$ and $\bbs\colon\bbd\to\bbk$ are
$\BO$-functors, the weighted limit $\{\bbf,\bbs\}$ in $\bbk$ exists
just when the presheaf
$[\bbd,\bbbo](\bbf,\bbk(-,\bbs))\colon\bbd\op\to\bbbo$ is representable.

A weighted colimit in $\bbk$ is a weighted limit in $\bbk\op$. As
usual, however, rather than consider a weight $\bbf\colon\bbd\to\bbbo$
and diagram $\bbs\colon\bbd\to\bbk\op$, we think of the diagram as
having the form $\bbs\colon\bbd\op\to\bbk$. In fact we more often (as
is also usual) replace $\bbd$ by its opposite, so that we have a presheaf $\bbf\colon\bbd\op\to\bbbo$ and diagram $\bbs\colon\bbd\to\bbk$.

In what follows, we work through what this means in various cases,
concentrating on the colimits, and leading up to the notion of
Kleisli object which appears in our main theorem.

\subsection{Tight colimits}

A tight colimit is a weighted colimit for which the weight is tight,
in the sense of Definition~\ref{defn:tight-weight}.

Consider a 2-category $\cd$, a presheaf $F\colon\cd\op\to\Cat$, and the corresponding
tight weight $F^{\id} \colon \id_*\cd\op\to\bbbo$ of Section~\ref{sect:CoBpresh1}.
By adjointness, a diagram $\bbs\colon\id_*\cd\to\bbk$ corresponds to a
2-functor $S\colon\cd\to\ck_t$ (its tight part); the loose part is
necessarily $E_\bbk S\colon\cd\to\ck_\ell$.

\begin{proposition}\label{prop:tight-colimit1}
  To give a $\BO$-enriched weighted colimit $F^{\id}*\bbs$ in the
  $\BO$-category $\bbk$ is equivalently to give a $\Cat$-enriched
  weighted colimit $F*S$ in $\ck_t$, with the further property that
  the induced functor
  \[ \ck_\ell(E_\bbk (F*S),B) \to [\cd\op,\Cat](F,\ck_\ell(E_\bbk
    S,B)) \] is fully faithful.
\end{proposition}

\begin{proof}
  Since
  $[\cd\op,\bbbo](F^{\id},\bbk(\bbs,B))_t\cong
  [\cd\op,\Cat](F,\ck(S,B))$, the tight aspect of the universal
  property says precisely that the $\BO$-colimit $F^{\id}*\bbs$ in
  $\bbk$ should be the $\Cat$-colimit $F*S$ in $\ck$. 

  To understand the loose part of the universal property, suppose that
  $F*S$ exists in $\ck_t$, and consider the following diagram.
  \begin{equation}\label{eq:loose-univ-prop} \xymatrix @C5pc {
      \ck_t(F*S,B) \ar[d]_{\cong} \ar[r]^{E_\bbk} & \ck_\ell(E_\bbk(F*S),B)
      \ar[d]^{\theta}  \\ 
      [\id_*\cd\op,\bbbo](F^{\id},\bbk(\bbs,B))_t \ar[d]_{\cong} \ar[r]^-{E_{[\id_*\cd\op,\bbbo]}} &
      [\id_*\cd\op,\bbbo](F^{\id},\bbk(\bbs,B))_\ell \ar[d]^{\psi}  \\
      [\cd\op,\Cat](F,\ck_t(S,B)) \ar[r]_-{[\cd\op,\Cat](F,E_\bbk)} &
      [\cd\op,\Cat](F,\ck_\ell(E_\bbk S,B)) }
  \end{equation}
  In the lower square, the left vertical is part of the change-of-base
  isomorphism of Section~\ref{sect:CoBpresh1}, while
  $E_{[\id_*\cd\op,\bbbo]}$ is the identity on objects and the right
  vertical is fully faithful by Section~\ref{sect:presheaf-cats}. In
  the upper square, the left vertical is the isomorphism expressing
  the universal property of the colimit $F*S$ (that is, the tight
  aspect of the universal property of $F^{\id}*\bbs$), while the loose
  aspect says that $\theta$ is invertible. Since the other three maps
  in the top square are bijective on objects, so is $\theta$. On the
  other hand, since $\psi$ is fully faithful, $\theta$ will be fully
  faithful (and so invertible) if and only if $\psi\theta$ is fully
  faithful.  
\end{proof}

\begin{proposition}
  In the case where $[\cd\op,\Cat](F,-)\colon[\cd\op,\Cat]\to\Cat$
  sends pointwise bijective-on-objects maps to bijective-on-objects
  maps, $F^{\id}*\bbs$ exists in $\bbk$ just when $F*S$ exists in
  $\ck_t$ and is preserved by $E\colon \ck_t\to\ck_\ell$.
\end{proposition}

\begin{proof}
  Suppose that $F*S$ exists in $\ck_t$. It will be preserved by
  $E_\bbk$ just when the displayed map in
  Proposition~\ref{prop:tight-colimit1} is invertible. This would be
  equivalent to it being fully faithful, as in the condition in
  Proposition~\ref{prop:tight-colimit1}, just when it is already known
  to be bijective on objects. But it is the image under
  $[\cd\op,\Cat](F,-)$ of the pointwise identity-on-objects map
  $E_\bbk\colon \ck_t(S,B)\to \ck_\ell(E_\bbk S,B)$, thus will indeed
  be bijective on objects under the hypotheses of the proposition.
\end{proof}

\begin{example}\label{ex:tight-colim-easy}
  Suppose that $\cd$ is a locally discrete 2-category --- it has no
  non-identity 2-cells --- and that $F=\Delta1\colon
  \cd\op\to\Cat$. Then
  $[\cd\op,\Cat](\Delta1,-)\colon[\cd\op,\Cat]\to\Cat$ does send
  pointwise bijective-on-object maps to bijective-on-object ones.

  Thus $\bbk$ has conical colimits indexed by ordinary categories if
  and only if $\ck_t$ has them and they are preserved by
  $\ck_t\to\ck_\ell$.

  Dually, $\bbk$ has conical limits indexed by ordinary categories if
  and only if $\ck_t$ has them and they are preserved by
  $\ck_t\to\ck_\ell$. 
\end{example}

\begin{example}
  On the other hand, consider the case of copowers by categories. Then
  $\cd=1$ and $F\colon 1\to\Cat$ picks out a category $X$. It is not
  the case that $\Cat(X,-)$ preserves bijective-on-object maps unless
  $X$ is discrete. A diagram $1\to\bbk$ corresponds to an object
  $A\in\ck$. The copower by $X$ will be a copower $X\cdot A$ in
  $\ck_t$ with the property that the induced
  \[ \ck_\ell(X\cdot A,B) \to \Cat(X,\ck_\ell(A,B)) \] is fully
  faithful.

  Thus we have bijections between maps $f\colon X\cdot A\to
  B$ in $\ck_t$ and functors $f'\colon X\to \ck_t(A,B)$; and
  furthermore, given $f,g\colon X\cdot A\to B$ and the corresponding
  $f',g'\colon X\to \ck_t(A,B)$ there are bijections between items in
  the left and right columns of the following table.
  
  \begin{tabular}{l|l}
    tight 2-cells $f\to g$ & natural transformations $f'\to g'$ \\ \hline
    loose 2-cells $f\to g$ & natural transformations $E_\bbk f'\to
                             E_\bbk g'$
  \end{tabular}
\end{example}

\subsection{Coproduct completions}

Our main result will concern free completions under a certain sort of
$\BO$-enriched Kleisli object. But perhaps it is also worth discussing
briefly the free completions under {\em coproducts}, seen as tight
colimits as in Example~\ref{ex:tight-colim-easy}.

Recall that the free completion of an ordinary category $\cx$ under
coproducts is given by $\Fam(\cx)$: an object is a ``family''
$X\colon I\to\cx$ of objects of $\cx$, where $I$ is an indexing set,
seen as a discrete category, and $X$ a functor; while a morphism has
the form
\[ \xymatrix{
    I \ar[dr]_X^{~}="1" \ar[rr]^-{f} && J \ar[dl]^Y_{~}="2" \\
    & \cx. \ar@{=>}"1";"2"^{F} } \]
If $\cx$ is not just a category but a 2-category, then there is a
natural way to make $\Fam(\cx)$ into a 2-category: a 2-cell
$(f,F)\to (g,G)$ can exist only if $f=g$, in which case it consists of
a modification $F\to G$. The resulting 2-category is the
free completion of $\cx$ under $\Cat$-enriched coproducts.

Since this
$\Cat$-enriched completion agrees with the ordinary one at the level
of underlying ordinary categories, if $E_\bba\colon\ca_t\to\ca_\ell$ acts as the
identity on objects and on 1-cells, then the same is true of
$\Fam(E)\colon \Fam(\ca_t)\to \Fam(\ca_\ell)$. This defines a
$\BO$-category $\Fam(\bba)$ for each $\BO$-category $\bba$.

By the universal property of the free completion, $\Fam(\ca_t)$ has
coproducts, and these are preserved by
$\Fam(E_\bba)\colon\Fam(\ca_t)\to \Fam(\ca_\ell)$. Thus the
$\BO$-category $\Fam(\bba)$ has tight coproducts. 

It follows easily from the universal property of the $\Cat$-enriched
$\Fam$ construction that $\Fam(\bba)$ is in fact the $\BO$-enriched
free completion under tight coproducts. 

Of course one can also modify this so as to deal with finite
coproducts (or indeed $\kappa$-small coproducts for some regular cardinal
$\kappa$) by limiting the size of the indexing sets~$I$.

\subsection{Loose colimits}

Once again, a loose colimit is a weighted colimit for which the weight
is loose, in the sense of Definition~\ref{defn:loose-weight}.

Consider a 2-presheaf $F\colon\cd\op\to\Cat$ and the corresponding
$\bbf\colon \disc_*\cd\op\to\bbbo$ as in Section~\ref{sect:CoBpresh2}
above.  By adjointness again, a diagram $\bbs\colon \disc_*\cd\to\bbk$
consists of just a 2-functor $S\colon \cd\to\ck_\ell$ --- we write
$S_t\colon \cd_1\to\ck_t$ for the uniquely determined tight part of $\bbs$, which
satisfies $E_\bbk S_t=S E_{\disc_*\cd}$ --- and a map
$\bbf\to \bbk(S,C)$ corresponds to a map $F\to \ck_\ell(S,C)$.

\begin{proposition}
  A $\BO$-enriched colimit $\bbf*\bbs$ in $\bbk$ is a $\Cat$-enriched
  colimit $F*S$ in $\ck_\ell$ with the property that, for maps
  $f,g\colon F*S\to B$ and a loose 2-cell $\phi\colon f\to g$, to give
  a tight structure to $\phi$ is equivalently to give a lifting
  $\phi'\colon f'\to g'$ as in 
  \[  \xymatrix @C3pc {
      F_1 \ar@{.>}@/^1pc/[rr]^-{f'}_(0.4){~}="1" \ar@{.>}@/_1pc/[rr]_-{g'}^(0.4){~}="2" \ar[d]_-{F_E} && \ck_t(S_t-,B) \ar[dd]^{E_\bbk} \\
      FE \ar[d] \\
      \ck_\ell(SE-,F*S) \ar@/^1pc/[rr]^-{\ck_\ell(SE-,f)}_(0.4){~}="3"
      \ar@/_1pc/[rr]_-{\ck_\ell(SE-,g)}^(0.4){~}="4" && \ck_\ell(SE-,B).
      \ar@{=>}"3";"4"^{\ck_\ell(SE-,\phi)}
      \ar@{:>}"1";"2"^{\phi'} 
    } \]


\end{proposition}

\begin{proof}
  First note that the existence and uniqueness of $f'$ and $g'$ are
  automatic, since $F_1$ is pointwise discrete and
  $E_\bbk\colon\ck_t(S_t-,B)\to \ck_\ell(SE-,B)$ is pointwise the identity on objects.
  
  By virtue of the isomorphism
  \[ [\disc_*\cd\op,\bbbo](\bbf,\bbk(\bbs,B))_\ell \cong
    [\cd\op,\Cat](F,\ck_\ell(S,B)) \] the loose aspect of the
  universal property says precisely that $\bbf*\bbs$ should be given
  by a $\Cat$-enriched limit $F*S$ in $\ck_\ell$. As for the tight aspect, since
  \[ \xymatrix{ [\disc_*\cd\op,\bbbo](\bbf,\bbk(\bbs,B))_t \ar[dd]
      \ar[r] &
      [\cd\op,\Cat](F,\ck_\ell(S,B)) \ar[d] \\
      &  [\cd\op_1,\Cat](FE,\ck_\ell (SE,B)) \ar[d]   \\
      [\cd\op_1,\Cat](F_1,\ck_t(S_t,B)) \ar[r] &
      [\cd\op_1,\Cat](F_1,\ck_\ell (SE,B)) } \] is a pullback as in Section~\ref{sect:presheaf-cats}, it
  follows that
  \[ \xymatrix{
      \ck_t(F*S,B) \ar[r] \ar[dd] & \ck_\ell(F*S,B) \ar[d]^-\cong \\
      & [\cd\op,\Cat](F,\ck_\ell(S,B)) \ar[d] \\
      [\cd\op_1,\Cat](F_1,\ck_t(S_t,B)) \ar[r] &
      [\cd\op_1,\Cat](F_1,\ck_\ell(SE,B)) } \] will need to be a
  pullback of categories. Since the two horizontals are both bijective
  on objects, this reduces to the condition in the proposition.
\end{proof}




\begin{remark}
  If $\cd$ is locally discrete, and $F$ takes values in discrete
  categories, then in fact $\bbf*\bbs$ is also a tight colimit, and just
  amounts to a colimit $F*S_t$ in $\ck_t$ which is preserved by
  $E_\bbk\colon\ck_t\to\ck_\ell$.
\end{remark}

\begin{remark}\label{rmk:loose-colimit-in-tight}
  If on the other hand $\bbk$ is tight, so that $\ck_t=\ck_\ell$ and
  $E_\bbk$ is the identity, then a loose colimit in $\bbk$ is just a
  (2-categorical) colimit in $\ck_\ell$. 
\end{remark}

\subsection{Kleisli objects}

Our main example involves the special case of the previous section
arising from 2-categorical Kleisli objects.
Let $\cm$ be the universal 2-category containing a monad, so that for
any 2-category $\ck$, there is a natural bijection between 2-functors
from $\cm$ to $\ck$ and monads in $\ck$. Now let $F\colon\cm\op\to\Cat$ be 
the weight for Kleisli objects, so that for a 2-functor
$T\colon\cm\to\ck$, a colimit $F*T$ is exactly a Kleisli object for the
monad corresponding to $T$. For details concerning $\cm$ and $F$, see
for example \cite[Section~5]{Street-Cat-limits} or
\cite[Section~8.2]{companion}.

Then $\cm_1$ is the universal category containing an endomorphism (in
other words, the one-object category whose morphisms are the natural
numbers with composition given by addition). And
$F_1\colon\cm\op_1\to\Cat$ picks out the discrete category $\bbz_{>0}$
of positive integers together with the successor endomorphism.

A diagram $\disc_*\cm\to\bbk$ consists of a monad $(A,t)$ in
$\ck_\ell$; we might call this a {\em loose monad in $\bbk$}. An
$\bbf$-weighted colimit of it is a Kleisli object $e\colon A\to A_t$
for the monad in $\ck_\ell$, with the additional property that if
$f,g\colon A_t\to B$ and $\phi\colon f\to g$ is loose, then
restriction along $e$ defines a bijection between liftings of $\phi$
to a tight 2-cell and liftings of $\phi e\colon ge\to he$ to a tight
2-cell.

\begin{remark}
  Of course this extra tightness condition is automatic if
  $\bbk=\id_*\ck$, so that all 2-cells are tight.
\end{remark}

\begin{proposition}\label{prop:enhanced-Kleisli}
  Let $e\colon A\to A'$ have a right adjoint $e\dashv r$ in
  $\ck_\ell$, and let $t$ be the induced monad in $\ck_\ell$; that is,
  the induced loose monad in $\bbk$. Then $e$ exhibits $A'$ as the
  Kleisli object in $\bbk$ if and only if
  \begin{enumerate}[(i)]
  \item $e$ exhibits $A'$ as the Kleisli object in $\ck_\ell$, and
  \item the square \[ \xymatrix{
        \ck_t(A_t,B) \ar[r]^{\ck_t(e,B)} \ar[d]_E & \ck_t(A,B) \ar[d]^E \\
        \ck_\ell(A_t,B) \ar[r]_{\ck_\ell(e,B)} & \ck_\ell(A,B) } \] is
    a pullback in $\Cat$, for each $B\in\bbk$.
  \end{enumerate}
\end{proposition}

\begin{proof}
  Here (i) expresses the loose aspect of the universal property. In
  the square in (ii), the vertical maps labelled $E$ are the identity
  on objects, so the object part of the pullback property is always
  true. What is left is the tight aspect of the universal property of
  the Kleisli object.
\end{proof}

\begin{definition}
  In this case we say that $e\dashv r$ is {\em an adjunction of enhanced
    Kleisli type in $\bbk$}. We call the colimit the {\em enhanced
    Kleisli object}, or just the Kleisli object if the context makes
  clear we are dealing with the $\BO$-enriched notion.
\end{definition}


Adjunctions of enhanced Kleisli type clearly compose, since Kleisli adjunctions compose in $\ck_\ell$ and
pullback squares can be pasted to give pullback squares.

\begin{proposition}\label{Kleisli-presheaf}
  Let $P=(P_t,P_\ell,P_E)$ be a presheaf on the $\BO$-category 
  $\bbk=(\ck_t,\ck_\ell,E_\bbk)$. A loose monad on $P$ consists of a
  monad $(S,m,i)$ on $P_\ell$ in $[\ck\op_\ell,\Cat]$, together with a
  lifting of the endomorphism $S\colon P_\ell\to P_\ell$ to a map
  $\overline{S}\colon P_t\to P_t$ in $[\ck\op_t,\Cat]$ making the diagram
  \[ \xymatrix{
      P_t \ar[r]^-{\overline{S}} \ar[d]_{P_E} & P_t \ar[d]^{P_E} \\
      P_\ell E_\bbk \ar[r]_-{S E_\bbk} & P_\ell E_\bbk } \] commute.

  The enhanced Kleisli object $P'=(P'_t,P'_\ell,P'_E)$ has $P'_t=P_t$,
  $P'_\ell$ the Kleisli object in $[\ck\op_\ell,\Cat]$ of $S$, and
  $P'_E= e_SE_\bbk.P_E$, where $e_S\colon P_\ell\to P'_\ell$ is the Kleisli
  map for $S$. The identity $1\colon P_t\to P'_t$ and
  $e_S\colon P_\ell\to P'_\ell$ define the Kleisli map $P\to P'$.
\end{proposition}

\begin{proof}
  The description of loose monads on $P$ follows immediately from the
  definitions. We derive the description of the enhanced Kleisli
  object using the characterization of these given in Proposition~\ref{prop:enhanced-Kleisli}.

  First we construct the loose adjunction $e\dashv r$. The right
  adjoint $r\colon P'\to P$ is defined by $S\colon P'_t=P_t\to P_t$ and
  $r_S\colon P'_\ell\to P_\ell$, where $r_S$ is the right adjoint of
  $e_S$. The unit and counit of $e\dashv r$ are just the unit and
  counit of the adjunction $e_S\dashv r_S$ in
  $[\ck\op_\ell,\Cat]$. This is clearly a loose adjunction and induces the
  original loose monad. We have to show that it is of enhanced Kleisli
  type.

  First suppose given a map $g\colon P\to Q$ with a loose opaction
  $gS\to g$; that is, an opaction $g_\ell S\to g_\ell$. Then $g_\ell$
  induces a unique $g'_\ell\colon P'_\ell\to Q_\ell$, which together
  with $g'_t:=g_t$ determines a map $g'\colon P'\to Q$. This defines a
  bijection between the objects of $[\bbk\op,\bbbo]_\ell(P',Q)$ and
  the objects of $[\bbk\op,\bbbo]_\ell(P,Q)$ equipped with an opaction
  of $S$, and this bijection extends to an isomorphism of categories
  exhibiting $P'$ as the Kleisli object in $[\bbk\op,\bbbo]_\ell$.

  It remains to check the tight aspect. This says that if
  $x,y\colon P'\to Q$ and $\xi\colon x_\ell\to y_\ell$, then to give a
  lifting of $\xi$ to some $\overline{\xi}\colon x_t\to y_t$ is
  equivalent to giving a lifting
  $\xi e_\ell\colon x_\ell e_\ell\to y_\ell e_\ell$ to some
  $x_t\to y_t$, which is trivially true.
\end{proof}



For a 2-category $\ck$, let $\KL(\ck)$ be the free completion of $\ck$
under ($\Cat$-enriched) Kleisli objects. As described in
\cite[Section~1]{ftm2}, this is equipped with a 2-functor
$\Mnd_*(\ck)\to\KL(\ck)$ which acts as the identity on objects and on
1-cells. Here $\Mnd_*(\ck)=\Mnd(\ck\op)\op$, and similarly
$\KL(\ck)=\EM(\ck\op)\op$.

Thus $\Mnd_*(\ck)\to\KL(\ck)$ can be seen as a $\BO$-category.

\begin{remark}
  In particular, if we take $\ck=\Cat$, then the $\BO$-category
  corresponding to $\Mnd_*(\Cat)\to\KL(\Cat)$ is isomorphic to the full
  sub-2-category of $\BO$ consisting of those identity-on-objects
  functors which are left adjoints. 
\end{remark}

\begin{theorem}
  The free completion of $\id_*\ck$ under enhanced Kleisli objects for loose
  monads is the $\BO$-category $\Mnd_*(\ck)\to\KL(\ck)$.
\end{theorem}

\begin{proof}
  The free completion $\widehat{\ck}$ will be given by the closure of the
  representables in $[\id_*\ck,\bbbo]$ under enhanced Kleisli
  objects.
  
  As observed above, the enhanced Kleisli-type adjunctions in a
  $\BO$-category are closed under composition, and so this closure
  under Kleisli objects will be of the ``one-step'' variety, and we
  can simply consider (the full subcategory consisting of) those
  objects which are Kleisli objects of loose monads on
  representables. Thus we can take the objects of $\widehat{\ck}$ to be
  the loose monads on representables in $[\id_*\ck,\bbbo]$, or
  equivalently the loose monads in $\id_*\ck$, or equivalently the
  monads in $\ck$.

  Given monads $(A,t)$ and $(B,s)$, the hom $\widehat{\ck}((A,t),(B,s))$
  will be the hom in $[\id_*\ck,\bbbo]$ between the Kleisli objects.

  By Proposition~\ref{Kleisli-presheaf} and the fact that
  $\ck_\ell=\ck_t$, the Kleisli object of $(A,t)$ is just the Kleisli
  object and Kleisli morphism
  $\ck_\ell(-,A) \to \ck_\ell(-,A)_{\ck_\ell(-,t)}$ in
  $[\ck\op_\ell,\Cat]$. We can now read off the various morphisms and
  2-cells. In particular, the loose 2-category $\widehat{\ck}_\ell$ is
  precisely $\KL(\ck)$, by exactly the calculation that was given in
  \cite{ftm2}. A morphism $(A,t)\to (B,s)$ is a pair $(f,\phi)$, where
  $f\colon A\to B$ and $\phi\colon ft\to sf$, subject to two
  equations. Given another morphism $(g,\psi)\colon (A,t)\to (B,s)$, a
  loose 2-cell $(f,\phi)\to (g,\psi)$ consists of a 2-cell
  $\rho\colon f\to sg$, subject to a single equation. To make this
  into a tight 2-cell is to give $\overline{\rho}\colon f\to g$ with
  $\eta g.\overline{\rho}=\rho$. But $\eta g.\overline{\rho}$ is a
  2-cell in $\KL(\ck)$ just when $\overline{\rho}\colon f\to g$
  defines a 2-cell $(f,\phi)\to (g,\psi)$ in $\Mnd_*(\ck)$.
\end{proof}

Dually, we have the notion of {\em enhanced Eilenberg-Moore object},
and the corresponding theorem.

\begin{theorem}\label{thm:EM}
  The free completion of $\id_*\ck$ under enhanced Eilenberg-Moore
  objects for loose monads is the $\BO$-category $\Mnd(\ck)\to\EM(\ck)$.
\end{theorem}

This is the promised universal property of
$\Mnd(\ck)\to\EM(\ck)$. Having fulfilled our promise, we conclude with
a few observations about Eilenberg-Moore objects as adjoints to the
inclusion $\ck\to\Mnd(\ck)$.

Since $\EM(\ck)$ is the free completion of $\ck$ under Eilenberg-Moore
objects, $\ck$ will have Eilenberg-Moore objects just when the
inclusion $\ck\to\EM(\ck)$ has a right adjoint. But why should
$\ck\to\Mnd(\ck)$ also have a right adjoint, as observed already in \cite{ftm}?

By Remark~\ref{rmk:loose-colimit-in-tight}, the tight $\BO$-category
$\id_*\ck$ has enhanced Eilenberg-Moore objects for loose monads if
and only if $\ck$ has ordinary Eilenberg-Moore objects for monads. By
Theorem~\ref{thm:EM}, this will be the case if and only if the
inclusion 
\[ \xymatrix{
    \ck \ar[r] \ar[d] & \Mnd(\ck) \ar[d] \\ \ck \ar[r] & \EM(\ck) } \]
has a right $\BO$-adjoint. Such an adjoint will imply in particular
that the horizontal components each have right $\Cat$-adjoints. Thus
we recover the result of \cite{ftm} that if $\ck$ has Eilenberg-Moore
objects then $\ck\to\Mnd(\ck)$ has a right adjoint.

\end{document}